\documentclass[11pt]{amsart}
\usepackage{amsmath,amssymb,mathrsfs,stmaryrd}
\newtheorem{theorem}{Theorem}
\newtheorem{proposition}[theorem]{Proposition}
\newtheorem{corollary}[theorem]{Corollary}
\newtheorem{lemma}[theorem]{Lemma}
\begin{document} 

\title[Isoperimetric and Weingarten surfaces]{Isoperimetric and Weingarten surfaces in the Schwarzschild manifold}
\author{Simon Brendle and Michael Eichmair}
\begin{abstract}
We show that any star-shaped convex hypersurface with constant Weingarten curvature in the deSitter-Schwarzschild manifold is a sphere of symmetry. Moreover, we study an isoperimetric problem for bounded domains in the doubled Schwarzschild manifold. We prove the existence of an isoperimetric surface for any value of the enclosed volume, and we completely describe the isoperimetric surfaces for very large enclosed volume. This complements work in H. Bray's thesis, where isoperimetric surfaces homologous to the horizon are studied.
\end{abstract}
\address{Department of Mathematics \\ Stanford University \\ Stanford, CA 94305 \\ U.S.A.}
\address{Departement Mathematik \\ ETH Z\"urich \\ 8092 Z\"urich \\ Switzerland}
\thanks{The first-named author was supported in part by the U.S. National Science Foundation under grant DMS-0905628. The second-named author was supported in part by the U.S. National Science Foundation under grant DMS-0906038 and by the Swiss National Science Foundation under grant SNF 200021-140467.}

\maketitle

\section{Introduction}

The classical Alexandrov theorem asserts that any closed embedded hypersurface in $\mathbb{R}^n$ with constant mean curvature is a round sphere. This theorem has been generalized by many authors. In particular, an analogue of Alexandrov's theorem holds in hyperbolic space (cf. \cite{Hijazi-Montiel-Roldan}, \cite{Montiel-Ros}), as well as in pseudo-hyperbolic space (see \cite{Montiel}). In a recent paper \cite{Brendle}, the first-named author proved a uniqueness theorem for constant mean curvature hypersurfaces in the deSitter-Schwarzschild manifold. Let us recall the definition of the deSitter-Schwarzschild manifold. Fix an integer $n \geq 3$, a real number $m >0$, and a real number $\kappa$ such that either $\kappa \leq 0$ or $\frac{n^n m^2 \kappa^{n-2}}{4 (n-2)^{n-2}} < 1$. Moreover, let 
\[I = \{s > 0 : 1 - m \, s^{2-n} - \kappa \, s^2 >0\}.\] 
Note that $I$ is a non-empty open interval, so we may write $I = (\underline{s},\overline{s})$. The deSitter-Schwarzschild manifold $(M,\bar{g})$ is defined by $M = S^{n-1} \times I$ and 
\[\bar{g} = \frac{1}{1 - m \, s^{2-n} - \kappa \, s^2} \, ds \otimes ds + s^2 \, g_{S^{n-1}}.\] 
Note that the metric extends smoothly to $S^{n-1} \times [\underline{s},\overline{s})$, and that the boundary component $S^{n-1} \times \{\underline{s}\}$ is totally geodesic. We will refer to this boundary component as the horizon of $(M, \bar{g})$.

\begin{theorem}[S.~Brendle \cite{Brendle}]
\label{alexandrov.theorem.in.schwarzschild}
Let $\Sigma$ be a closed, embedded hypersurface in the deSitter-Schwarzschild manifold $(M,\bar{g})$. If $\Sigma$ has constant mean curvature, then $\Sigma$ is a slice $S^{n-1} \times \{s\}$. 
\end{theorem}

We note that surfaces of constant mean curvature play an important role in general relativity; see e.g. \cite{Bray-thesis}, \cite{Christodoulou-Yau}, \cite{Huisken-Yau}, \cite{Eichmair-Metzger1}, \cite{Eichmair-Metzger2}, \cite{Eichmair-Metzger3}, \cite{Qing-Tian}.

It is interesting to replace the mean curvature of $\Sigma$ by other functions of the principal curvatures. In this direction, Ros \cite{Ros} showed that any closed, embedded hypersurface in $\mathbb{R}^n$ with constant $\sigma_p$ is a round sphere. Here, $\sigma_p$ denotes the $p$-th elementary symmetric polynomial in the principal curvatures. This result was generalized to hyperbolic space by Montiel and Ros \cite{Montiel-Ros}. Finally, He, Li, Ma, and Ge \cite{He-Li-Ma-Ge} studied the case of anisotropic higher order mean curvatures. 

We first analyze surfaces in deSitter-Schwarzschild space with constant higher order mean curvature. Under some extra assumptions, we are able to show that such surfaces are spheres of symmetry:  

\begin{theorem}
\label{higher.order.mean.curvatures}
Let $\Sigma$ be a closed, embedded hypersurface in the deSitter-Schwarzschild manifold $(M,\bar{g})$ that is star-shaped and convex. Moreover, suppose that $\sigma_p = \text{\rm constant}$, where $\sigma_p$ denotes the $p$-th elementary symmetric polynomial in the principal curvatures. Then $\Sigma$ is a slice $S^{n-1} \times \{s\}$.
\end{theorem}

The convexity assumption is needed to control certain curvature terms arising in the Codazzi equations. Theorem \ref{higher.order.mean.curvatures} is a special case of a stronger result that applies to more general warped product manifolds. We will explain this in Section \ref{weingarten.surfaces}.

We now consider the case when $\kappa = 0$. In this case, $\overline{s} = \infty$ and $(M, \bar g)$  is the standard Schwarzschild manifold. By reflection across the totally geodesic boundary component $S^{n-1} \times \{\underline{s}\}$, we obtain a complete manifold $(\bar{M},\bar{g})$ which, up to scaling, is isometric to $\mathbb{R}^n \setminus \{0\}$ equipped with the metric $\bar{g}_{ij} = (1 + |x|^{2-n})^{\frac{4}{n-2}} \, \delta_{ij}$. In this representation, the horizon is the coordinate sphere $\partial B_1(0)$. We will refer to $(\bar{M},\bar{g})$ as the doubled Schwarzschild manifold. In his thesis, Bray \cite{Bray-thesis} studied isoperimetric surfaces in the Schwarzschild manifold which are homologous to the horizon (see also \cite{Bray-Morgan}):

\begin{theorem}[H.~Bray \cite{Bray-thesis}]
\label{bray.theorem}
Let $\Sigma$ be a sphere of symmetry in the doubled Schwarzschild manifold. Then $\Sigma$ has least area among all comparison surfaces that are homologous to $\Sigma$ and which enclose the same oriented volume with the horizon.
\end{theorem}

In view of this result, it is natural to wonder about isoperimetric surfaces in the doubled Schwarzschild manifold that are null-homologous. Our first result here establishes that isoperimetric regions exist in $(\bar{M}, \bar{g})$ for every volume:

\begin{theorem}
\label{isoperimetric.1}
Given any $V>0$, there exists a bounded Borel set $\Omega$ of finite perimeter and volume $V$ in the doubled Schwarzschild manifold that has least perimeter amongst all such sets.
\end{theorem}

Classical results in geometric measure theory show that the reduced boundary $\partial^*\Omega$ of a set $\Omega$ as in Theorem \ref{isoperimetric.1} is a smooth volume preserving stable constant mean curvature hypersurface, such that $\partial^*\Omega$ is relatively open in $\overline {\partial^*\Omega}$, and such that the Hausdorff dimension of $\overline {\partial^*\Omega} \setminus \partial^*\Omega$ does not exceed $n-8$. Moreover, one can (and we always will) choose a representative of $\Omega$ such that $\partial \Omega = \overline {\partial^*\Omega}$. \\

In the 1980's, S.T. Yau asked whether there exist constant mean curvature surfaces in the doubled Schwarzschild manifold other than the spheres of symmetry. By choosing $V>0$ very small in Theorem \ref{isoperimetric.1} we obtain examples of such surfaces that are even isoperimetric. We observe that the existence of small surfaces of constant mean curvature in the doubled Schwarzschild manifold can alternatively be deduced from general perturbation results of Pacard and Xu \cite{Pacard-Xu}. The construction in \cite{Pacard-Xu} neither implies nor indicates that their surfaces are isoperimetric. 

Finally, we give a precise description of the large volume isoperimetric regions in the doubled Schwarzschild manifold. A crucial ingredient in the proof is an effective version of Theorem \ref{isoperimetric.1} established in \cite{Eichmair-Metzger1} (when $n=3$) and \cite{Eichmair-Metzger3} (for $n \geq 3$). 

\begin{theorem}
\label{isoperimetric.2}
Let $\Omega$ be an isoperimetric region in the doubled Schwarzschild manifold. If the volume of $\Omega$ is sufficiently large, then $\Omega$ is bounded by two spheres of symmetry.
\end{theorem}

In fact, there are exactly two isoperimetric regions for every given large volume, and they are obtained from each other by reflection across the horizon. Under the extra assumption that $\Omega$ has smooth boundary, we can prove the following:

\begin{theorem} 
\label{isoperimetric.3}
Suppose that $\Omega$ is a smooth isoperimetric region in the doubled Schwarzschild manifold. Then either $\Omega$ is bounded by two spheres of symmetry, or the boundary of $\Omega$ is connected and intersects the horizon.
\end{theorem}

If $3 \leq n < 8$, the smoothness assumption in Theorem \ref{isoperimetric.3} is always satisfied. We expect that Theorem \ref{alexandrov.theorem.in.schwarzschild} can be generalized to constant mean curvature surfaces with a small singular set, so that the smoothness assumption in Theorem \ref{isoperimetric.3} can be dropped. 

We note that similar results for the cylinder have been obtained by Pedrosa \cite{Pedrosa} using symmetrization techniques. 

The authors would like to thank Professors Jan Metzger, Sebastian Montiel, Brian White, and Shing-Tung Yau for discussions and their interest. 

\section{Weingarten surfaces in warped product manifolds}

\label{weingarten.surfaces}

Let us fix an integer $n \geq 3$. We consider the manifold $M = S^{n-1} \times [0,\bar{r})$ equipped with a Riemannian metric of the form $\bar{g} = dr \otimes dr + h(r)^2 \, g_{S^{n-1}}$. We assume that the warping function $h: [0,\bar{r}) \to \mathbb{R}$ satisfies the following conditions: 
\begin{itemize} 
\item[(H1)] $h'(0) = 0$ and $h''(0) > 0$.
\item[(H2)] $h'(r) > 0$ for all $r \in (0,\bar{r})$.
\item[(H3)] The function 
\[2 \, \frac{h''(r)}{h(r)} - (n-2) \, \frac{1 - h'(r)^2}{h(r)^2}\] 
is non-decreasing for $r \in (0,\bar{r})$.
\item[(H4)] We have $\frac{h''(r)}{h(r)} + \frac{1-h'(r)^2}{h(r)^2} > 0$ for all $r \in (0,\bar{r})$.
\end{itemize}
We note that the Ricci and scalar curvature of $(M,\bar{g})$ are given by 
\begin{align} 
\label{ricci.tensor}
\text{\rm Ric} 
&= -\Big ( \frac{h''(r)}{h(r)} - (n-2) \, \frac{1-h'(r)^2}{h(r)^2} \Big ) \, \bar{g} \notag \\ 
&- (n-2) \, \Big ( \frac{h''(r)}{h(r)} + \frac{1-h'(r)^2}{h(r)^2} \Big ) \, dr \otimes dr 
\end{align} 
and 
\begin{equation} 
R = -(n-1) \, \Big ( 2 \, \frac{h''(r)}{h(r)} - (n-2) \, \frac{1 - h'(r)^2}{h(r)^2} \Big ). 
\end{equation}
Hence, condition (H3) is equivalent to saying that the scalar curvature is a non-increasing function of $r$. Moreover, condition (H4) is equivalent to saying that the Ricci curvature is smallest in the radial direction. In particular, the conditions (H1)--(H4) are all satisfied on the deSitter-Schwarzschild manifolds.

Note that a closed hypersurface $\Sigma$ in $M$ is either null-homologous or  bounds a compact region together with $S^{n-1} \times \{0\}$. We say that $\Sigma$ is star-shaped if there is a choice of unit normal $\nu$ such that $\langle \frac{\partial}{\partial r},\nu \rangle \geq 0$ on $\Sigma$. We say that $\Sigma$ is convex if there is a choice of unit normal such that all its principal curvatures are non-negative. 

The following is the main result of this section: 

\begin{theorem}
\label{weingarten}
Let $(M,\bar{g})$ be a warped product manifold satisfying conditions (H1)--(H4) above. Let $\Sigma$ be a closed, embedded hypersurface in $(M,\bar{g})$ that is star-shaped and convex. Moreover, suppose that $\sigma_p = \text{\rm constant}$, where $\sigma_p$ denotes the $p$-th elementary symmetric polynomial in the principal curvatures. Then $\Sigma$ is a slice $S^{n-1} \times \{r\}$ for some $r \in (0,\bar{r})$.
\end{theorem}

As in \cite{Brendle}, we define a function $f$ and a vector field $X$ by $f = h'(r)$ and $X = h(r) \, \frac{\partial}{\partial r}$. Note that $X$ is a conformal vector field; in fact, $\bar{D} X = f \, \bar{g}$.

We now consider a hypersurface $\Sigma$ in $M$. Let $\{e_1,\hdots,e_{n-1}\}$ be a local orthonormal frame on $\Sigma$, and let $\nu$ denote the unit normal to $\Sigma$. Moreover, let $h_{ij} = \langle \bar{D}_{e_i} \nu,e_j \rangle$ denote the second fundamental form of $\Sigma$, and let $\sigma_p$ denote the $p$-th elementary symmetric polynomial in the eigenvalues of $h$. Finally, we put $T_{ij}^{(p)} = \frac{\partial}{\partial h_{ij}} \sigma_p$. We may view $T_{ij}^{(p)}$ as a symmetric two-tensor on $\Sigma$. It turns out that the divergence of $T_{ij}^{(p)}$ has a special structure (see also \cite{Viaclovsky}):

\begin{proposition}
\label{preparation}
Suppose that $\Sigma$ star-shaped and convex. Then 
\[\sum_{i=1}^{n-1} \langle X,e_i \rangle \, (D_{e_j} T^{(p)})(e_i,e_j) \geq 0.\] 
Here, $D$ denotes the Levi-Civita connection on $\Sigma$.
\end{proposition}

\textbf{Proof.} 
We may write 
\[\sum_{p=0}^{n-1} t^p \, \sigma_p = \det(I+th).\] 
Differentiating this identity with respect to $h_{ij}$, we obtain 
\[\sum_{p=1}^{n-1} t^p \, T_{ij}^{(p)} = t \, \det(I+th) \, G(t)_{ij},\] 
where $G(t)$ denotes the inverse of $I+th$. We now take the divergence on both sides of this identity. This yields 
\begin{align*} 
\sum_{j, p=1}^{n-1} t^p \, D_j T_{ij}^{(p)} 
&= t \, \det(I+th) \, \sum_{j=1}^{n-1} D_j G(t)_{ij} \\ 
&+ t^2 \, \det(I+th) \, \sum_{j=1}^{n-1} G(t)_{ij} \, G(t)_{kl} \, D_j h_{kl} \\ 
&= -t^2 \, \det(I+th) \, \sum_{j,k,l=1}^{n-1} G(t)_{ik} \, G(t)_{jl} \, D_j h_{kl} \\ 
&+ t^2 \, \det(I+th) \, \sum_{j=1}^{n-1} G(t)_{ij} \, G(t)_{kl} \, D_j h_{kl} \\ 
&= -t^2 \, \det(I+th) \, \sum_{j,k,l=1}^{n-1} G(t)_{ik} \, G(t)_{jl} \, (D_j h_{kl} - D_k h_{jl}). 
\end{align*} 
Using the Codazzi equations, we obtain 
\[D_j h_{kl} - D_k h_{jl} = R(e_j,e_k,e_l,\nu),\] 
where $R$ denotes the Riemann curvature tensor of $(M,\bar{g})$. Since $\bar{g}$ is locally conformally flat, the curvature tensor of $\bar{g}$ is given by $\frac{1}{n-2} \, A \owedge \bar{g}$, where $A$ is the Schouten tensor of $\bar{g}$. Therefore, 
\[D_j h_{kl} - D_k h_{jl} = -\frac{1}{n-2} \, (\text{\rm Ric}(e_j,\nu) \, \bar{g}(e_k,e_l) - \text{\rm Ric}(e_k,\nu) \, \bar{g}(e_j,e_l)).\] 
Putting these facts together, we obtain 
\begin{align*} 
&(n-2) \sum_{j, p=1}^{n-1}  t^p \, D_j T_{ij}^{(p)} \\ 
&= t^2 \, \det(I+th) \, \sum_{j,k=1}^{n-1} G(t)_{ik} \, G(t)_{jk} \, \text{\rm Ric}(e_j,\nu) \\ 
&- t^2 \, \det(I+th) \, \text{\rm tr}(G(t)) \, \sum_{j=1}^{n-1} G(t)_{ij} \, \text{\rm Ric}(e_j,\nu), 
\end{align*} 
and hence 
\begin{align*} 
&(n-2) \sum_{i,j, p=1}^{n-1} t^p \, \langle X,e_i \rangle \, D_j T_{ij}^{(p)} \\ 
&= t^2 \, \det(I+th) \, \sum_{i,j,k=1}^{n-1} G(t)_{ik} \, G(t)_{jk} \, \langle X,e_i \rangle \, \text{\rm Ric}(e_j,\nu) \\ 
&- t^2 \, \det(I+th) \, \text{\rm tr}(G(t)) \, \sum_{i,j=1}^{n-1} G(t)_{ij} \, \langle X,e_i \rangle \, \text{\rm Ric}(e_j,\nu). 
\end{align*} 
Without loss of generality, we may assume that $h_{ij}$ is diagonal with eigenvalues $\lambda_1,\hdots,\lambda_{n-1} \geq 0$. Then 
\begin{align*} 
&(n-2) \sum_{i,j, p=1}^{n-1} t^p \, \langle X,e_i \rangle \, D_j T_{ij}^{(p)} \\ 
&= -t^2 \, \det(I+th) \, \sum_{i \neq j} \frac{1}{(1+t\lambda_i)(1+t\lambda_j)} \, \langle X,e_j \rangle \, \text{\rm Ric}(e_j,\nu) \\ 
&= -t^2 \, \sum_{i \neq j} \bigg ( \prod_{k \in \{1,\hdots,n-1\} \setminus \{i,j\}} (1+t\lambda_k) \bigg ) \, \langle X,e_j \rangle \, \text{\rm Ric}(e_j,\nu) \\ 
&= -\sum_{p=2}^{n-1} \sum_{j=1}^{n-1} (n-p) \, t^p \, \sigma_{p-2}(\lambda_1,\hdots,\lambda_{j-1},\lambda_{j+1},\hdots,\lambda_{n-1}) \, \langle X,e_j \rangle \, \text{\rm Ric}(e_j,\nu).
\end{align*} 
Comparing coefficients gives 
\begin{align*} 
&\sum_{i,j=1}^{n-1} \langle X,e_i \rangle \, D_j T_{ij}^{(p)} \\ 
&= -\frac{n-p}{n-2} \, \sum_{j=1}^{n-1} \sigma_{p-2}(\lambda_1,\hdots,\lambda_{j-1},\lambda_{j+1},\hdots,\lambda_{n-1}) \, \langle X,e_j \rangle \, \text{\rm Ric}(e_j,\nu) 
\end{align*} 
for each $p \in \{2,\hdots,n-1\}$. On the other hand, it follows from (\ref{ricci.tensor}) and (H4) that $\text{\rm Ric}(e_j,\nu)$ is a negative multiple of $\langle X,e_j \rangle \, \langle X,\nu \rangle$. Since $\langle X,\nu \rangle \geq 0$, we conclude that $\langle X,e_j \rangle \, \text{\rm Ric}(e_j,\nu) \leq 0$ for $j = 1, \hdots, n-1$. From this, the assertion follows. \\

The following result can be viewed as an analogue of the Minkowski-type formula established in \cite{Brendle}, see also \cite{Alias-deLira-Malacarne}, Section 8:

\begin{proposition}
\label{minkowski}
Suppose that $\Sigma$ is star-shaped and convex. Then 
\[p \int_\Sigma \langle X,\nu \rangle \, \sigma_p \geq (n-p) \int_\Sigma f \, \sigma_{p-1}.\] 
\end{proposition}

\textbf{Proof.} 
Let $\xi$ denote the orthogonal projection of $X$ to the tangent space of $\Sigma$, i.e. 
\[\xi = X - \langle X,\nu \rangle \, \nu.\] 
Then 
\[D_i \xi_j = \bar{D}_i X_j -  \langle X,\nu \rangle \, h_{ij}  = f \, g_{ij} -  \langle X,\nu \rangle \, h_{ij}.\] 
Hence 
\[\sum_{i,j=1}^{n-1} D_i(\xi_j \, T_{ij}^{(p)}) = f \sum_{i=1}^{n-1} T_{ii}^{(p)} - \sum_{i,j=1}^{n-1} T_{ij}^{(p)} \langle X,\nu \rangle \, h_{ij}  + \sum_{i,j=1}^{n-1} \xi_j \, D_i T_{ij}^{(p)}.\] 
Since $\sigma_p$ is a homogeneous function of degree $p$, we have 
\[\sum_{i,j=1}^{n-1} h_{ij} \, T_{ij}^{(p)} = p \, \sigma_p\] 
by Euler's theorem. Moreover, it is easy to see that 
\[\sum_{i=1}^{n-1} T_{ii}^{(p)} = (n-p) \, \sigma_{p-1}.\] 
Finally, we have 
\[\sum_{i,j=1}^{n-1} \xi_j \, D_i T_{ij}^{(p)} \geq 0\] 
by Proposition \ref{preparation}. Putting these facts together, we obtain 
\[\sum_{i,j=1}^{n-1} D_i(\xi_j \, T_{ij}^{(p)}) \geq (n-p) \, f \, \sigma_{p-1} - p \,  \langle X,\nu \rangle \, \sigma_p.\] 
Hence, the assertion follows from the divergence theorem. \\

After these preparations, we are now able to complete the proof of Theorem \ref{weingarten}. Suppose that $\Sigma$ is a  star-shaped and convex hypersurface with the property that $\sigma_p$ is constant. By Proposition \ref{minkowski}, we have 
\[p \int_\Sigma  \langle X,\nu \rangle \, \sigma_p \geq (n-p) \int_\Sigma f \, \sigma_{p-1}.\] 
Since $\sigma_p$ is constant, it follows that 
\[p \int_\Sigma \langle X,\nu \rangle \geq (n-p) \int_\Sigma f \, \frac{\sigma_{p-1}}{\sigma_p}.\] 
Using the Newton inequality, we obtain 
\[(n-p) \, \sigma_{p-1} \, \sigma_1 \geq (n-1)p \, \sigma_p.\] 
Therefore, 
\begin{equation} 
\label{ineq.1}
\int_\Sigma \langle X,\nu \rangle \geq (n-1) \int_\Sigma \frac{f}{H}, 
\end{equation}
where $H = \sigma_1$ denotes the mean curvature of $\Sigma$. On the other hand, it was shown in \cite{Brendle}, Section 3, that 
\begin{equation} 
\label{ineq.2}
(n-1) \int_\Sigma \frac{f}{H} \geq \int_\Sigma \langle X,\nu \rangle.  
\end{equation}
Therefore, equality holds in (\ref{ineq.1}) and (\ref{ineq.2}). Since equality holds in (\ref{ineq.2}), it follows from results in \cite{Brendle} that $\Sigma$ is a slice $S^{n-1} \times \{r\}$ for some $r \in (0,\bar{r})$.

\section{Null-homologous isoperimetric surfaces in the doubled Schwarzschild manifold}

In this section, we consider the doubled Schwarzschild manifold $(\bar{M}, \bar{g}) = (\mathbb{R}^n \setminus \{0\}, (1+|x|^{2-n})^{\frac{4}{n-2}} \, \delta_{ij})$ discussed in the introduction. 

Given any $V>0$, we define $A_{\bar{g}}(V)$ as the infimum of $\mathscr{H}_{\bar{g}}^{n-1}(\partial^* \Omega)$ where $\Omega$ ranges over all Borel subsets of $(\bar{M},\bar{g})$ with finite perimeter and $\text{\rm vol}_{\bar{g}}(\Omega) = V$. If such a set $\Omega$ realizes the infimum, i.e. if $A_{\bar{g}}(\text{\rm vol}_{\bar{g}}(\Omega)) = \mathscr{H}_{\bar{g}}^{n-1}(\partial^* \Omega)$, then $\Omega$ is called an isoperimetric region.

\begin{lemma}
\label{boundedness}
If $\Omega$ is an isoperimetric region, then $\Omega$ is bounded.
\end{lemma} 

\textbf{Proof.} 
Suppose that $\Omega$ is unbounded. Then we can find a sequence of points $p_k$ in the support of $\Omega$ such that $\text{dist}_{\bar{g}}(p_k,p_l) > 2$ for $k \neq l$. Let $B_k$ denote the geodesic ball of radius $k^{-\frac{1}{n-1}}$ centered at $p_k$. It follows from the monotonicity formula that $\liminf_{k \to \infty} k \, \mathscr{H}_{\bar{g}}^{n-1}( B_k \cap \partial^* \Omega) > 0$. This implies that $\mathscr{H}_{\bar{g}}^{n-1}(\partial^* \Omega) \geq \sum_{k=1}^\infty \mathscr{H}_{\bar{g}}^{n-1}(B_k \cap \partial^* \Omega) = \infty$, a contradiction. \\

The behavior of minimizing sequences for the isoperimetric problem in general is described in \cite{Ritore-Rosales}, Theorem 2.1. In conjunction with the characterization of isoperimetric regions in Euclidean space, the following result was obtained in \cite{Eichmair-Metzger1}, Proposition 4.2:

\begin{proposition} \label{minimizing.sequence} 
Given any $V > 0$, there exists a (possibly empty) isoperimetric region $\Omega \subset \bar{M}$ and a real number $\rho \geq 0$ such that 
\[\text{\rm vol}_{\bar{g}}(\Omega) + \frac{1}{n} \, \omega_{n-1} \, \rho^n = V\] 
and 
\[\mathscr{H}_{\bar{g}}^{n-1}(\partial^* \Omega) + \omega_{n-1} \, \rho^{n-1} = A(V).\] 
\end{proposition}

We now establish the existence of isoperimetric regions in $(\bar{M}, \bar{g})$ of any given volume. We note that in \cite{Ritore}, M. Ritor\'e has constructed examples of complete rotationally symmetric Riemannian surfaces in which no solutions of the isoperimetric problem exist for any volume.  We first establish an auxiliary results:

\begin{proposition} 
\label{profile}
Given $V>0$ and any compact set $K \subset \bar{M}$, we can find a smooth region $D \subset \bar{M} \setminus K$ such that $\text{\rm vol}_{\bar{g}}(D) = V$ and $\mathscr{H}_{\bar{g}}^{n-1}(\partial D) < (n^{n-1} \, \omega_{n-1})^{\frac{1}{n}} \, V^{\frac{n-1}{n}}$.
\end{proposition}

\textbf{Proof.} 
Let us consider a ball $B = \{x: |x-a| \leq r\}$, where $r$ is bounded and $|a|$ is large compared to $r$. By Corollary \ref{isoperimetric.ratio}, we can perturb $B$ in a suitable way to obtain a region $D$ with 
\begin{align*} 
\mathscr{H}_{\bar{g}}^{n-1}(\partial D) 
&= (n^{n-1} \, \omega_{n-1})^{\frac{1}{n}} \, \text{\rm vol}_{\bar{g}}(D)^{\frac{n-1}{n}} \\ 
&\cdot \bigg [ 1 
- \frac{2(n-2)(n-1)^2}{(n+1)(n+2)(n+4)} \, \frac{r^4}{|a|^{2n}} + O(|a|^{-2n-1}) \bigg] .
\end{align*}
In particular, we have that
\[\mathscr{H}_{\bar{g}}^{n-1}(\partial D) < (n^{n-1} \, \omega_{n-1})^{\frac{1}{n}} \, \text{\rm vol}_{\bar{g}}(D)^{\frac{n-1}{n}}\] 
provided $|a|$ is sufficiently large. The construction of $D$ in Appendix \ref{offcenter} is continuous in $r$, and the volume of $D$ converges to $\frac{1}{n} \, \omega_{n-1} \, r^n$ as $|a| \to \infty$. Hence, given any $a$ sufficiently large, we can choose $r$ such that $\text{\rm vol}_{\bar{g}}(D) = V$. \\

\begin{theorem} 
\label{existence.of.isoperimetric.surfaces} 
For every $V > 0$, the infimum $A(V)$ is achieved. 
\end{theorem} 

\textbf{Proof.} 
By Proposition \ref{minimizing.sequence}, we can find an isoperimetric region $\Omega \subset \bar{M}$ and a real number $\rho \geq 0$ such that 
\[\text{\rm vol}_{\bar{g}}(\Omega) + \frac{1}{n} \, \omega_{n-1} \, \rho^n = V\] 
and 
\[\mathscr{H}_{\bar{g}}^{n-1}(\partial^* \Omega) + \omega_{n-1} \, \rho^{n-1} = A(V).\] 
By Lemma \ref{boundedness}, $\Omega$ is bounded. 

We claim that $\rho=0$. Indeed, if $\rho>0$, Proposition \ref{profile} implies the existence of a smooth bounded region $D$ whose closure is disjoint from the closure of $\Omega$, and which satisfies $\text{\rm vol}_{\bar{g}}(D) = \frac{1}{n} \, \omega_{n-1} \, \rho^n$ and $\mathscr{H}_{\bar{g}}^{n-1}(\partial D) < \omega_{n-1} \, \rho^{n-1}$. Consequently, we have 
\[\text{\rm vol}_{\bar{g}}(\Omega \cup D) = \text{\rm vol}_{\bar{g}}(\Omega) + \frac{1}{n} \, \omega_{n-1} \, \rho^n = V\] 
and 
\[\mathscr{H}_{\bar{g}}^{n-1}(\partial^* \Omega \cup \partial D) < \mathscr{H}_{\bar{g}}^{n-1}(\partial^*\Omega) + \omega_{n-1} \, \rho^{n-1} = A(V).\] 
This contradicts the definition of $A(V)$. Hence $\rho = 0$ and $\Omega$ is an isoperimetric region of volume $V$. \\

We next describe the isoperimetric regions with large volume. We need the following effective version of a result of H. Bray's (Theorem \ref{bray.theorem}) from \cite{Eichmair-Metzger3}: 

\begin{proposition} [\cite{Eichmair-Metzger3}, Proposition 3.4] 
\label{effective} 
Given $(\tau,\eta) \in (1,\infty) \times (0,1)$ there exists $V_0 > 0$ so that the following holds: Let $\Omega$ be a bounded Borel set with finite perimeter in the doubled Schwarzschild manifold $(\bar{M},\bar{g}) = (\mathbb{R}^n \setminus \{0\}, (1 +|x|^{2-n})^{\frac{4}{n-2}} \, \delta_{ij})$, and let $r \geq 1$ be such that 
\[\text{\rm vol}_{\bar{g}}(\Omega \setminus B_1(0)) = \text{\rm vol}_{\bar{g}}(B_r(0) \setminus B_1(0)) \geq V_0.\] 
If $\Omega$ is $(\tau,\eta)$-off-center, i.e. if $\mathscr{H}_{\bar{g}}^{n-1}(\partial^* \Omega \setminus B_{\tau r}(0)) \geq \eta \, \mathscr{H}_{\bar{g}}^{n-1}(\partial B_r(0))$, then 
\begin{equation} 
\label{comparison.in.Schwarzschild}
\mathscr{H}_{\bar{g}}^{n-1}(\partial^* \Omega \setminus B_1(0)) \geq \mathscr{H}_{\bar{g}}^{n-1}(\partial B_r(0)) + c \eta \Big ( 1 - \frac{1}{\tau} \Big )^2 \, r.
\end{equation}
Here, $c >0$ is a constant that only depends on $n$. 
\end{proposition}

\begin{theorem}
\label{large.isoperimetric.regions}
There exists $V_1>0$ with the following property: If $\Omega$ is an isoperimetric region in the doubled Schwarzschild manifold $(\bar{M},\bar{g})$ with $\text{\rm vol}_{\bar{g}}(\Omega) \geq V_1$, then $\partial^*\Omega$ is a union of two spheres of symmetry.
\end{theorem}

\textbf{Proof.} 
Suppose this is false. Let $\Omega_k$ be a sequence of isoperimetric regions with $\text{\rm vol}_{\bar{g}}(\Omega_k) \to \infty$ such that $\partial^*\Omega_k$ is not a union of two spheres of symmetry. Reflecting across the horizon if necessary, we may assume  that $\text{\rm vol}_{\bar{g}}(\Omega_k \setminus B_1(0)) \to \infty$. Let $r_k \geq 1$ be such that $\text{\rm vol}_{\bar{g}}(\Omega_k \setminus B_1(0)) = \text{\rm vol}_{\bar{g}}(B_{r_k}(0) \setminus B_1(0))$. Since $\Omega_k$ is an isoperimetric region, it follows that 
\begin{equation} 
\label{inequality}
\mathscr{H}_{\bar{g}}^{n-1}(\partial^*\Omega_k \setminus B_1(0)) \leq \mathscr{H}_{\bar{g}}^{n-1}(\partial B_{r_k}(0)) + \mathscr{H}_{\bar{g}}^{n-1}(\partial B_1(0)). 
\end{equation} 
We now consider two cases: 

\textit{Case 1:} Suppose first that 
\[\liminf_{k \to \infty} r_k^{-n} \, \text{\rm vol}_{\bar{g}}(\Omega_k \setminus B_{\tau r_k}(0)) = 0\] 
for every $\tau > 1$. As in the proof of Theorem 5.1 in \cite{Eichmair-Metzger1} we see that, possibly after passing to a subsequence, $\Omega_k \subset B_{2r_k}(0)$ and that the rescaled regions $r_k^{-1} \, \Omega_k$ converge to the unit ball $B_1(0)$ in Euclidean space away from the origin. In particular, we have 
\[B_{r_k/2}(0) \setminus B_{r_k/4}(0) \subset \Omega_k \subset B_{2r_k}(0)\] 
for some large integer $k$. By Allard's theorem, $\overline {\partial^* \Omega_k} \setminus B_{r_k/2}(0)$ is a smooth connected constant mean curvature surface and hence, by Theorem \ref{alexandrov.theorem.in.schwarzschild}, a centered coordinate sphere. Let $\hat{r}_k = \inf \{r \in (0,r_k/2): B_{r_k/2}(0) \setminus B_r(0) \subset \Omega_k\}$. By the half-space theorem, the coordinate sphere $\partial B_{\hat{r}_k}(0)$ intersects $\overline{\partial^* \Omega_k}$ in the regular set $\partial^* \Omega_k$. The maximum principle then implies that $\hat{r}_k < 1$. Theorem \ref{bray.theorem} gives that $\overline {\partial^*\Omega_k} \cap B_1(0)$ is a centered coordinate sphere. Therefore, $\Omega_k$ is smooth and its boundary is a union of two coordinate spheres. This contradicts the choice of $\Omega_k$. 

\textit{Case 2:} We now assume that there exists a real number $\tau > 1$ such that 
\[\liminf_{k \to \infty} r_k^{-n} \, \text{\rm vol}_{\bar{g}}(\Omega_k \setminus B_{\tau r_k}(0)) > 0.\] 
This implies that 
\[\liminf_{k \to \infty} r_k^{1-n} \, \mathscr{H}_{\bar{g}}^{n-1}(\partial^*\Omega_k \setminus B_{\tau r_k}(0)) > 0\] 
for some number $\tau > 1$. Consequently, we can find a real number $\eta \in (0,1)$ with the property that the sets $\Omega_k$ are $(\tau,\eta)$-off-center when $k$ is sufficiently large. Using Proposition \ref{effective}, we obtain 
\[\mathscr{H}_{\bar{g}}^{n-1}(\partial^*\Omega_k \setminus B_1(0)) \geq \mathscr{H}_{\bar{g}}^{n-1}(\partial B_{r_k}(0)) + c \eta \Big ( 1 - \frac{1}{\tau} \Big )^2 \, r_k\] 
for $k$ sufficiently large. This contradicts (\ref{inequality}). \\

\begin{corollary}
Let $\Omega$ be an isoperimetric region in the doubled Schwarzschild manifold. If the volume of $\Omega$ is sufficiently large, then $\Omega = B_{r_1}(0) \setminus B_{r_0}(0)$ for suitable real numbers $r_1 > 1 > r_0$. Moreover, $r_0r_1 \neq 1$ and 
\[\frac{r_1 \, (r_1^{n-2}-1)}{(r_1^{n-2}+1)^{\frac{n}{n-2}}} = \frac{r_0 \, (1-r_0^{n-2})}{(1+r_0^{n-2})^{\frac{n}{n-2}}}.\] 
\end{corollary}

\textbf{Proof.} 
By Theorem \ref{large.isoperimetric.regions}, the boundary $\partial^*\Omega$ is a union of  two spheres of symmetry. Since the components of $\partial^*\Omega$ have the same positive mean curvature, we conclude that $\Omega = B_{r_1}(0) \setminus B_{r_0}(0)$ where $r_1 > 1 > r_0$. Moreover, we have 
\[\frac{r_1 \, (r_1^{n-2}-1)}{(r_1^{n-2}+1)^{\frac{n}{n-2}}} = \frac{r_0 \, (1-r_0^{n-2})}{(1+r_0^{n-2})^{\frac{n}{n-2}}} = \frac{H}{n-1},\] 
where $H$ denotes the mean curvature of $\partial^*\Omega$. It remains to show that $r_0r_1 \neq 1$. Indeed, if $r_0r_1 = 1$, then $r_1$ is large, and we have $\text{\rm vol}_{\bar{g}}(\Omega) = \frac{2}{n} \, \omega_{n-1} \, r_1^n + O(r_1^{n-1})$ and $\mathscr{H}_{\bar{g}}^{n-1}(\partial^*\Omega) = 2 \, \omega_{n-1} \, r_1^{n-1} + O(r_1^{n-2})$. On the other hand, by Corollary \ref{isoperimetric.ratio}, we have that $A(V) < (n^{n-1} \, \omega_{n-1})^{\frac{1}{n}} \, V^{\frac{n-1}{n}}$. This gives a contradiction. 
\\

Finally, we describe the proof of Theorem \ref{isoperimetric.3}. Let $\Omega$ be an isoperimetric region in the doubled Schwarzschild manifold $(\bar{M},\bar{g})$ with smooth boundary. We consider two cases: 

\textit{Case 1:} Suppose that one of the components of $\partial \Omega$ does not intersect the horizon. By Theorem \ref{alexandrov.theorem.in.schwarzschild}, this component must be a centered coordinate sphere. A maximum principle argument as in the proof of Theorem \ref{large.isoperimetric.regions} shows $\Omega$ is bounded by two spheres of symmetry.

\textit{Case 2:} Suppose next that every component of $\partial \Omega$ intersects the horizon. If $\Omega$ is disconnected, we may take one connected component of $\Omega$ and rotate it until it touches another connected component of $\Omega$. This process does not change the isoperimetric property of the region since volume and boundary area stay unchanged. Clearly, the final configuration is not optimal for the isoperimetric problem. This is a contradiction. Therefore, $\Omega$ must be connected. Moreover, if the complement $\bar{M} \setminus \Omega$ has two unbounded components $D_1$ and $D_2$, then the maximum principle implies that the boundary $\partial D_1$ must lie on one side of the horizon, contrary to our assumption. Consequently, $\bar{M} \setminus \Omega$ has exactly one unbounded component. Elementary topological consideration now imply that the boundary $\partial \Omega$ is connected. This completes the proof of Theorem \ref{isoperimetric.3}. 

\appendix 

\section{The isoperimetric ratio of coordinate balls in the doubled Schwarzschild manifold} 

\label{offcenter}

Let us consider the Riemannian metric 
\[\bar{g}_{ij} = (1+|y+a|^{2-n})^{\frac{4}{n-2}} \, \delta_{ij}\] 
on $\mathbb{R}^3 \setminus \{-a\}$. Let $B_r$ denote the coordinate ball of radius $r$ centered at the origin. We want to analyze the isoperimetric ratio of $B_r$ with respect to the metric $\bar{g}$ when $r$ is bounded and $|a|$ is large compared to $r$.

\begin{proposition}
\label{coordinate.balls}
We have 
\begin{align*} 
\mathscr{H}_{\bar{g}}^{n-1}(\partial B_r) 
&= \omega_{n-1} \, r^{n-1} \, (1+|a|^{2-n})^{\frac{2(n-1)}{n-2}} \\ 
&\cdot \bigg [ 1 + \frac{(n-1) \, r^2}{|a|^{2n-2}} \, (1 + |a|^{2-n})^{-2} + \frac{n(n-1)^2}{2(n+2)} \, \frac{r^4}{|a|^{2n}} + O(|a|^{-2n-1}) \bigg ] 
\end{align*}
and 
\begin{align*} 
\text{\rm vol}_{\bar{g}}(B_r) 
&= \frac{1}{n} \, \omega_{n-1} \, r^n \, (1+|a|^{2-n})^{\frac{2n}{n-2}} \\ 
&\cdot \bigg [ 1 + \frac{n \, r^2}{|a|^{2n-2}} \, (1+|a|^{2-n})^{-2} + \frac{n^2(n-1)}{2(n+4)} \, \frac{r^4}{|a|^{2n}} + O(|a|^{-2n-1}) \bigg ] 
\end{align*}
as $|a| \to \infty$.
\end{proposition}

\textbf{Proof.} 
It follows from Taylor's theorem that 
\begin{align} 
\label{taylor.a}
&(1 + |y+a|^{2-n})^{\frac{2(n-1)}{n-2}} \notag \\ 
&= (1 + |a|^{2-n})^{\frac{2(n-1)}{n-2}} \notag \\ 
&+ \frac{2(n-1)}{n-2} \, (1 + |a|^{2-n})^{\frac{n}{n-2}} \, (|a+y|^{2-n} - |a|^{2-n}) \notag \\ 
&+ \frac{n(n-1)}{(n-2)^2} \, (1 + |a|^{2-n})^{\frac{2}{n-2}} \, (|a+y|^{2-n} - |a|^{2-n})^2 \\ 
&+ \frac{2n(n-1)}{3 \, (n-2)^3} \, (1 + |a|^{2-n})^{-\frac{n-4}{n-2}} \, (|a+y|^{2-n} - |a|^{2-n})^3 \notag \\ 
&+ O \big ( (|a+y|^{2-n} - |a|^{2-n})^4 \big ) \notag
\end{align} 
and 
\begin{align} 
\label{taylor.b}
&(1 + |y+a|^{2-n})^{\frac{2n}{n-2}} \notag \\ 
&= (1 + |a|^{2-n})^{\frac{2n}{n-2}} \notag \\ 
&+ \frac{2n}{n-2} \, (1 + |a|^{2-n})^{\frac{n+2}{n-2}} \, (|a+y|^{2-n} - |a|^{2-n}) \notag \\ 
&+ \frac{n(n+2)}{(n-2)^2} \, (1 + |a|^{2-n})^{\frac{4}{n-2}} \, (|a+y|^{2-n} - |a|^{2-n})^2 \\ 
&+ \frac{4n(n+2)}{3 \, (n-2)^3} \, (1 + |a|^{2-n})^{-\frac{n-6}{n-2}} \, (|a+y|^{2-n} - |a|^{2-n})^3 \notag \\ 
&+ O \big ( (|a+y|^{2-n} - |a|^{2-n})^4 \big ). \notag
\end{align} 
The mean value property of harmonic functions implies that 
\begin{equation} 
\label{integral.1.a}
\int_{\partial B_r} (|y+a|^{2-n} - |a|^{2-n}) = 0 
\end{equation}
and 
\begin{equation} 
\label{integral.1.b}
\int_{B_r} (|y+a|^{2-n} - |a|^{2-n}) = 0. 
\end{equation}
We next observe that 
\begin{align*} 
&|y+a|^{2-n} - |a|^{2-n} \\ 
&= -(n-2) \, \frac{\langle a,y \rangle}{|a|^n} - \frac{n-2}{2} \, \frac{|a|^2 \, |y|^2 - n \, \langle a,y \rangle^2}{|a|^{n+2}} \\ 
&+ \frac{n(n-2)}{6} \, \frac{3 \, |a|^2 \, |y|^2 \, \langle a,y \rangle - (n+2) \, \langle a,y \rangle^3}{|a|^{n+4}} + O(|a|^{-n-2}). 
\end{align*} 
This implies 
\begin{align} 
\label{integral.2.a}
&\int_{\partial B_r} (|y+a|^{2-n} - |a|^{2-n})^2 \notag \\ 
&= (n-2)^2 \int_{\partial B_r} \frac{\langle a,y \rangle^2}{|a|^{2n}} + \frac{(n-2)^2}{4} \int_{\partial B_r} \frac{(|a|^2 \, |y|^2 - n \, \langle a,y \rangle^2)^2}{|a|^{2n+4}} \notag \\ 
&- \frac{n(n-2)^2}{3} \int_{\partial B_r} \frac{3 \, |a|^2 \, |y|^2 \, \langle a,y \rangle^2 - (n+2) \, \langle a,y \rangle^4}{|a|^{2n+4}} + O(|a|^{-2n-1}) \\ 
&= \frac{(n-2)^2}{n} \, \omega_{n-1} \, \frac{r^{n+1}}{|a|^{2n-2}} + \frac{(n-1)(n-2)^2}{2(n+2)} \, \omega_{n-1} \, \frac{r^{n+3}}{|a|^{2n}} + O(|a|^{-2n-1}) \notag 
\end{align} 
and 
\begin{align} 
\label{integral.2.b}
&\int_{B_r} (|y+a|^{2-n} - |a|^{2-n})^2 \notag \\ 
&= (n-2)^2 \int_{B_r} \frac{\langle a,y \rangle^2}{|a|^{2n}} + \frac{(n-2)^2}{4} \int_{B_r} \frac{(|a|^2 \, |y|^2 - n \, \langle a,y \rangle^2)^2}{|a|^{2n+4}} \notag \\ 
&- \frac{n(n-2)^2}{3} \int_{B_r} \frac{3 \, |a|^2 \, |y|^2 \, \langle a,y \rangle^2 - (n+2) \, \langle a,y \rangle^4}{|a|^{2n+4}} + O(|a|^{-2n-1}) \\ 
&= \frac{(n-2)^2}{n(n+2)} \, \omega_{n-1} \, \frac{r^{n+2}}{|a|^{2n-2}} + \frac{(n-1)(n-2)^2}{2(n+2)(n+4)} \, \omega_{n-1} \, \frac{r^{n+4}}{|a|^{2n}} + O(|a|^{-2n-1}). \notag
\end{align} 
Moreover, we have 
\begin{equation} 
\label{integral.3.a}
\int_{\partial B_r} (|a+y|^{2-n} - |a|^{2-n})^3 = O(|a|^{2-3n}) 
\end{equation}
and 
\begin{equation} 
\label{integral.3.b}
\int_{B_r} (|a+y|^{2-n} - |a|^{2-n})^3 = O(|a|^{2-3n}). 
\end{equation}
Using (\ref{taylor.a}), (\ref{integral.1.a}), (\ref{integral.2.a}), and (\ref{integral.3.a}), we obtain 
\begin{align*} 
\mathscr{H}_{\bar{g}}^{n-1}(\partial B_r) 
&= \int_{\partial B_r} (1 + |y+a|^{2-n})^{\frac{2(n-1)}{n-2}} \\ 
&= \omega_{n-1} \, r^{n-1} \, (1 + |a|^{2-n})^{\frac{2(n-1)}{n-2}} \\ 
&+ \frac{n(n-1)}{(n-2)^2} \, (1 + |a|^{2-n})^{\frac{2}{n-2}} \int_{\partial B_r} (|a+y|^{2-n} - |a|^{2-n})^2 \\ 
&+ O(|a|^{2-3n}) \\ 
&= \omega_{n-1} \, r^{n-1} \, (1 + |a|^{2-n})^{\frac{2(n-1)}{n-2}} \\ 
&+ (n-1) \, \omega_{n-1} \, (1 + |a|^{2-n})^{\frac{2}{n-2}} \, \frac{r^{n+1}}{|a|^{2n-2}} \\ 
&+ \frac{n(n-1)^2}{2(n+2)} \, \omega_{n-1} \, (1 + |a|^{2-n})^{\frac{2}{n-2}} \, \frac{r^{n+3}}{|a|^{2n}} + O(|a|^{-2n-1}). 
\end{align*} 
Similarly, it follows from (\ref{taylor.b}), (\ref{integral.1.b}), (\ref{integral.2.b}), and (\ref{integral.3.b}) that 
\begin{align*} 
\text{\rm vol}_{\bar{g}}(B_r) 
&= \int_{B_r} (1 + |y+a|^{2-n})^{\frac{2n}{n-2}} \\ 
&= \frac{1}{n} \, \omega_{n-1} \, r^n \, (1 + |a|^{2-n})^{\frac{2n}{n-2}} \\ 
&+ \frac{n(n+2)}{(n-2)^2} \, (1 + |a|^{2-n})^{\frac{4}{n-2}} \int_{B_r} (|a+y|^{2-n} - |a|^{2-n})^2 \\ 
&+ O(|a|^{2-3n}) \\ 
&= \frac{1}{n} \, \omega_{n-1} \, r^n \, (1 + |a|^{2-n})^{\frac{2n}{n-2}} \\ 
&+ (1 + |a|^{2-n})^{\frac{4}{n-2}} \, \omega_{n-1} \, \frac{r^{n+2}}{|a|^{2n-2}} \\ 
&+ \frac{n(n-1)}{2(n+4)} \, \omega_{n-1} \, (1 + |a|^{2-n})^{\frac{4}{n-2}} \, \frac{r^{n+4}}{|a|^{2n}} + O(|a|^{-2n-1}). 
\end{align*} 
This completes the proof. \\

\begin{lemma}
The mean curvature of $\partial B_r$ is given by 
\[H = \frac{n-1}{r} \, \bigg [ (1+|a|^{2-n})^{-\frac{2}{n-2}} - \frac{ |a|^2 \, |y|^2 - n \, \langle a,y \rangle^2}{|a|^{n+2}} \bigg ] + O(|a|^{-n-1}).\] 
\end{lemma} 

\textbf{Proof.} 
The standard formula for the change of the mean curvature under a conformal change of the metric gives 
\[H = \frac{n-1}{r} \, (1+|y+a|^{2-n})^{-\frac{2}{n-2}} \, \bigg [ 1 + \frac{2}{n-2} \sum_{i=1}^n y_i \, \frac{\partial}{\partial y_i} \log (1+|y+a|^{2-n}) \bigg ].\] 
Note that  
\begin{align*} 
&(1+|y+a|^{2-n})^{-\frac{2}{n-2}} \\ 
&= (1+|a|^{2-n})^{-\frac{2}{n-2}} \, \bigg [ 1 + (1+|a|^{2-n})^{-1} \, \Big ( 2 \, \frac{\langle a,y \rangle}{|a|^n} + \frac{|a|^2 \, |y|^2 - n \, \langle a,y \rangle^2}{|a|^{n+2}} \Big ) \bigg ] \\ 
&+ O(|a|^{-n-1}). 
\end{align*}
Similarly, 
\begin{align*} 
&\log(1+|y+a|^{2-n}) \\ 
&= \log(1 + |a|^{2-n}) - (n-2) \, (1+|a|^{2-n})^{-1} \, \Big ( \frac{\langle a,y \rangle}{|a|^n} + \frac{|a|^2 \, |y|^2 - n \, \langle a,y \rangle^2}{2 \, |a|^{n+2}} \Big ) \\ 
&+ O(|a|^{-n-1}).  
\end{align*}
Hence 
\begin{align*} 
&\sum_{i=1}^n y_i \, \frac{\partial}{\partial y_i} \log(1+|y+a|^{2-n}) \\ 
&= -(n-2) \, (1+|a|^{2-n})^{-1} \, \Big ( \frac{\langle a,y \rangle}{|a|^n} + \frac{|a|^2 \, |y|^2 - n \, \langle a,y \rangle^2}{|a|^{n+2}} \Big ) \\ 
&+ O(|a|^{-n-1}). 
\end{align*}
Putting these facts together, we obtain 
\begin{align*} 
&(1+|y+a|^{2-n})^{-\frac{2}{n-2}} \, \bigg [ 1 + \frac{2}{n-2} \sum_{i=1}^n y_i \, \frac{\partial}{\partial y_i} \log(1+|y+a|^{2-n}) \bigg ] \\ 
&= (1+|a|^{2-n})^{-\frac{2}{n-2}} \, \bigg [ 1 - (1+|a|^{2-n})^{-1} \, \frac{|a|^2 \, |y|^2 - n \, \langle a,y \rangle^2}{|a|^{n+2}} \bigg ] + O(|a|^{-n-1}). 
\end{align*} From this, the assertion follows. \\

Proposition \ref{coordinate.balls} shows that the isoperimetric ratio of the ball $B_r \subset (\bar{M},\bar{g})$ is greater than the isoperimetric ratio of a ball in Euclidean space. We overcome this obstacle by perturbing the coordinate ball $B_r$ in a suitable way. Let us define a function $f: \partial B_r \to \mathbb{R}$ by 
\[f(y) = \frac{(n-1) \, r}{n+1} \, \frac{|a|^2 \, r^2 - n \, \langle a,y \rangle^2}{|a|^{n+2}} + c,\] 
where the constant $c$ is chosen such that $\int_{\partial B_r} f \, d\mu_{\bar{g}} = 0$. Clearly, $c = O(|a|^{-n-1})$ and $f = O(|a|^{-n})$. We now consider the graph 
\[\Sigma = \{\exp_y(f(y) \, \nu(y)): y \in \partial B_r\},\] 
where $\exp$ denotes the exponential map with respect to $\bar{g}$ and $\nu$ denotes the unit normal to $\partial B_r$ with respect to $\bar{g}$. Moreover, let $\Omega$ denote the region enclosed by $\Sigma$.

\begin{proposition}
We have 
\begin{align*} 
\mathscr{H}_{\bar{g}}^{n-1}(\Sigma) 
&= \omega_{n-1} \, r^{n-1} \, (1+|a|^{2-n})^{\frac{2(n-1)}{n-2}} \\ 
&\cdot \bigg [ 1 + \frac{(n-1) \, r^2}{|a|^{2n-2}} \, (1 + |a|^{2-n})^{-2} \\ 
&\hspace{5mm} + \frac{n(n-1)^2(3n^2-6n+7)}{2(n+2)(n+1)^2} \, \frac{r^4}{|a|^{2n}} + O(|a|^{-2n-1}) \bigg ] 
\end{align*}
and 
\begin{align*} 
\text{\rm vol}_{\bar{g}}(\Omega) 
&= \frac{1}{n} \, \omega_{n-1} \, r^n \, (1+|a|^{2-n})^{\frac{2n}{n-2}} \\ 
&\cdot \bigg [ 1 + \frac{n \, r^2}{|a|^{2n-2}} \, (1+|a|^{2-n})^{-2} \\ 
&\hspace{5mm} + \frac{n(n-1)(3n^4+6n^3-13n^2+24n-8)}{2(n+2)(n+4)(n+1)^2} \, \frac{r^4}{|a|^{2n}} + O(|a|^{-2n-1}) \bigg ] 
\end{align*}
as $|a| \to \infty$.
\end{proposition} 

\textbf{Proof.} The surface area of $\Sigma$ is given by 
\begin{align*} 
\mathscr{H}_{\bar{g}}^{n-1}(\Sigma) 
&= \mathscr{H}_{\bar{g}}^{n-1}(\partial B_r) + \int_{\partial B_r} H \, f \, d\mu_{\bar{g}} \\ 
&+ \frac{1}{2} \int_{\partial B_r} \big ( |\nabla f|^2 + H^2 \, f^2 - |I\!I|^2 \, f^2 - \text{\rm Ric}(\nu,\nu) \, f^2 \big ) \, d\mu_{\bar{g}} + O(|a|^{-2n-1}). 
\end{align*}
Moreover, the volume of $\Omega$ satisfies 
\[\text{\rm vol}_{\bar{g}}(\Omega) = \text{\rm vol}_{\bar{g}}(B_r) + \frac{1}{2} \int_{\partial B_r} H \, f^2 \, d\mu_{\bar{g}} + O(|a|^{-2n-1}).\] 
Using the identity 
\[H = \frac{n-1}{r} \, \bigg [ (1+|a|^{2-n})^{-\frac{2}{n-2}} - \frac{ |a|^2 \, |y|^2 - n \, \langle a,y \rangle^2}{|a|^{n+2}} \bigg ] + O(|a|^{-n-1})\] 
and the relation $\int_{\partial B_r} f \, d\mu_{\bar{g}} = 0$, we obtain 
\begin{align*} 
\int_{\partial B_r} H \, f \, d\mu_{\bar{g}} 
&= -\frac{n-1}{r} \, \int_{\partial B_r} \frac{|a|^2 \, r^2 - n \, \langle a,y \rangle^2}{|a|^{n+2}} \, f \, d\mu_{\bar{g}} + O(|a|^{-2n-1}) \\ 
&= -\frac{(n-1)^2}{n+1} \int_{\partial B_r} \frac{(|a|^2 \, r^2 - n \, \langle a,y \rangle^2)^2}{|a|^{2n+4}} \, d\mu_{\bar{g}} + O(|a|^{-2n-1}). 
\end{align*} 
Moreover, the function $f$ satisfies 
\[\Delta_{\partial B_r} f = -\frac{2n(n-1)}{(n+1)r} \, \frac{|a|^2 \, r^2 - n \, \langle a,y \rangle^2}{|a|^{n+2}} + O(|a|^{-n-1}),\] 
hence 
\begin{align*} 
&\Delta_{\partial B_r} f - \frac{(n-2)(n-1)}{r^2} \, f \\ 
&= -\frac{(n^2-n+2)(n-1)}{(n+1)r} \, \frac{|a|^2 \, r^2 - n \, \langle a,y \rangle^2}{|a|^{n+2}} + O(|a|^{-n-1}). 
\end{align*}
Therefore, we have 
\begin{align*} 
&\int_{\partial B_r} \big ( |\nabla f|^2 + H^2 \, f^2 - |I\!I|^2 \, f^2 - \text{\rm Ric}(\nu,\nu) \, f^2 \big ) \, d\mu_{\bar{g}} \\ 
&= \int_{\partial B_r} \Big ( |\nabla f|^2 + \frac{(n-2)(n-1)}{r^2} \, f^2 \Big ) \, d\mu_{\bar{g}} + O(|a|^{-2n-1}) \\ 
&= -\int_{\partial B_r} \Big ( \Delta_{\partial B_r} f - \frac{(n-2)(n-1)}{r^2} \, f \Big ) \, f \, d\mu_{\bar{g}} + O(|a|^{-2n-1}) \\ 
&= \frac{(n^2-n+2)(n-1)^2}{(n+1)^2} \int_{\partial B_r} \frac{(|a|^2 \, r^2 - n \, \langle a,y \rangle^2)^2}{|a|^{2n+4}} \, d\mu_{\bar{g}} + O(|a|^{-2n-1}). 
\end{align*}
Putting these facts together, we obtain 
\begin{align*} 
&\mathscr{H}_{\bar{g}}^{n-1}(\Sigma) \\ 
&= \mathscr{H}_{\bar{g}}^{n-1}(\partial B_r) + \int_{\partial B_r} H \, f \, d\mu_{\bar{g}} \\ 
&+ \frac{1}{2} \int_{\partial B_r} \big ( |\nabla f|^2 + H^2 \, f^2 - |I\!I|^2 \, f^2 - \text{\rm Ric}(\nu,\nu) \, f^2 \big ) \, d\mu_{\bar{g}} + O(|a|^{-2n-1}) \\ 
&= \mathscr{H}_{\bar{g}}^{n-1}(\partial B_r) + \frac{n(n-3)(n-1)^2}{2(n+1)^2} \int_{\partial B_r} \frac{(|a|^2 \, r^2 - n \, \langle a,y \rangle^2)^2}{|a|^{2n+4}} \, d\mu_{\bar{g}} + O(|a|^{-2n-1}) \\
&= \mathscr{H}_{\bar{g}}^{n-1}(\partial B_r) + \frac{n(n-3)(n-1)^3}{(n+2)(n+1)^2} \, \omega_{n-1} \, \frac{r^{n+3}}{|a|^{2n}} + O(|a|^{-2n-1})  
\end{align*} 
and 
\begin{align*} 
\text{\rm vol}_{\bar{g}}(\Omega) 
&= \text{\rm vol}_{\bar{g}}(B_r) + \frac{n-1}{2r} \int_{\partial B_r} f^2 \, d\mu_{\bar{g}} + O(|a|^{-2n-1}) \\ 
&= \text{\rm vol}_{\bar{g}}(B_r) + \frac{(n-1)^3 \, r}{2(n+1)^2} \int_{\partial B_r} \frac{(|a|^2 \, r^2 - n \, \langle a,y \rangle^2)^2}{|a|^{2n+4}} \, d\mu_{\bar{g}} + O(|a|^{-2n-1}) \\ 
&= \text{\rm vol}_{\bar{g}}(B_r) + \frac{(n-1)^4}{(n+2)(n+1)^2} \, \omega_{n-1} \, \frac{r^{n+4}}{|a|^{2n}} + O(|a|^{-2n-1}). 
\end{align*}
Hence, the assertion follows from Proposition \ref{coordinate.balls}. \\

\begin{corollary}
\label{isoperimetric.ratio}
We have 
\begin{align*} 
\mathscr{H}_{\bar{g}}^{n-1}(\Sigma) 
&= (n^{n-1} \, \omega_{n-1})^{\frac{1}{n}} \, \text{\rm vol}_{\bar{g}}(\Omega)^{\frac{n-1}{n}} \\ 
&\cdot \bigg [ 1 
- \frac{2(n-2)(n-1)^2}{(n+1)(n+2)(n+4)} \, \frac{r^4}{|a|^{2n}} + O(|a|^{-2n-1}) \bigg ] 
\end{align*}
as $|a| \to \infty$.
\end{corollary}

\end{document}